\documentclass[12pt,reqno]{amsart}
\usepackage{amssymb,amscd,amsbsy}
\newcommand{\lb}{\linebreak}

\renewcommand{\a}{\alpha}

\newcommand{\z}{\zeta}

\newcommand{\vt}{\vartheta}

\newcommand{\s}{\sigma}

\newcommand{\f}{\varphi}

\renewcommand{\L}{\Lambda}

\renewcommand{\O}{\Omega}
\newcommand{\U}{\Upsilon}

\newcommand{\B}{{\mathcal B}}

\newcommand{\h}{{\mathcal H}}

\newcommand{\K}{{\mathcal K}}
\newcommand{\cL}{{\mathcal L}}
\newcommand{\M}{{\mathcal M}}
\newcommand{\N}{{\mathcal N}}

\newcommand{\V}{{\mathcal V}}
\newcommand{\W}{{\mathcal W}}

\newcommand{\C}{{\Bbb C}}
\newcommand{\T}{{\Bbb T}}
\newcommand{\pp}{{\Bbb P}}
\newcommand{\dd}{{\Bbb D}}

\newcommand{\Z}{{\Bbb Z}}
\newcommand{\mm}{{\Bbb M}}
\newcommand{\0}{{\boldsymbol{0}}}

\newcommand{\bs}{\boldsymbol}

\newcommand{\rf}[1]{(\ref{#1})}

\newcommand{\df}{\stackrel{\mathrm{def}}{=}}
\newcommand{\dist}{\operatorname{dist}}
\newcommand{\Ker}{\operatorname{Ker}}

\newcommand{\spn}{\operatorname{span}}

\newcommand{\rank}{\operatorname{rank}}

\newcommand{\eeq}{\end{equation}}
\newcommand{\beq}{\begin{equation}}
\newcommand{\bay}{\begin{eqnarray}}
\newcommand{\ba}{\begin{align*}}
\newcommand{\ea}{\end{align*}}
\newcommand{\ey}{\end{eqnarray}}
\newcommand{\bey}{\begin{eqnarray*}}
\newcommand{\eey}{\end{eqnarray*}}

\newcommand{\be}{\infty}

\newcommand{\bl}{\blacksquare}
\newcommand{\ess}{\operatorname{ess}}
\newcommand{\ind}{\operatorname{ind}}
\newcommand{\wind}{\operatorname{wind}}
\newcommand{\Range}{\operatorname{Range}}

\newcommand{\Pf}{{\bf Proof. }}

\newcommand{\ov}{\overline}

\newtheorem{thm}{\hspace{\parindent}Theorem}[section]

\theoremstyle{remark}

\newtheorem*{rem*}{Remark}

\newcommand{\co}{{\mathcal O}}
\newcommand\up{\upsilon}

\begin{document}

\newcommand{\vse}{\vspace{.2in}}
\numberwithin{equation}{section}

\title{\bf Approximation by analytic operator functions. 
Factorizations and very badly approximable functions}
\author{V.V. Peller and S.R. Treil}
\maketitle

\begin{abstract}
This is a continuation of our earlier paper \cite{PT3}. We consider here 
operator-valued functions (or infinite matrix functions)
on the unit circle $\T$ and study the problem of approximation
by bounded analytic operator functions. We discuss thematic and canonical factorizations
of operator functions and study badly approximable and very badly approximable operator
functions. 

We obtain algebraic and geometric characterizations of badly approximable and very badly
approximable operator functions. Note that there is an important difference between
the case of finite matrix functions and the case of operator functions. Our criteria for a
function to be very badly approximable in the case of finite
matrix functions also guarantee that the zero function is the only
superoptimal approximant. However in the case of operator functions this is not true.
\end{abstract}

\section{\bf Introduction}
\setcounter{equation}{0}

\

Our previous paper \cite{PT3} was devoted to a characterization of very badly approximable
matrix functions. In this paper we consider the case of operator-valued functions or,
which is equivalent, infinite matrix functions.

{\bf Background (best approximation and badly approximable functions).}
The classical problem of analytic approximation is for a given bounded function $\f$ 
on the unit circle $\T$ is to find a function $f$ in the Hardy class $H^\be$ 
such that
$$
\|\f - f\|_\infty = \dist_{L^\infty}(\f, H^\infty)= 
\inf_{h\in H^\infty} \|\f-h\|_\infty. 
$$
Such a best approximant $f$ always exist (a compactness argument) and as was proved by
S. Khavinson \cite{Kh} it is unique if $\f$ is continuous. 

A function $\f\in L^\infty$ is called {\em badly approximable} if 
$$
\|\f\|_\infty = \dist_{L^\infty}(\f, H^\infty),
$$
There is an elegant characterization of the set of continuous 
badly approximable functions:
{\it a nonzero continuous function $\f$ on $\T$
is badly approximable if and only if it has constant modulus and
its winding number $\wind\f$ with respect to the origin
is negative} (see \cite{AAK}, \cite {Po}). 

To extend this criterion to a broader class of functions $\f$, we need the 
notion of Hankel and Toeplitz operators. The {\it Toeplitz operator}
$T_\f:H^2\to H^2$ and the {\it Hankel operator} 
$H_\f:H^2\to H^2_-\df L^2\ominus H^2$ are defined by
\bay
\label{th}
T_\f f=\pp_+\f f,\quad H_\f f=\pp_-\f f,
\ey
where $\pp_-$ and $\pp_+$ are the orthogonal projections onto the subspaces $H^2$ and 
\lb$H^2_-\df L^2\ominus H^2$ of $L^2$.

It is well known (see e.g., \cite{D} or \cite{P}) 
that if $\f\in C(\T)$ and $\f$ does not vanish on $\T$, then 
the Toeplitz operator $T_\f$ on the Hardy class $H^2$ is Fredholm and \lb$\ind T_\f
=-\wind\f$ (recall that for a Fredholm operator $A$, its {\it index} is defined as
$\ind A=\dim\Ker A -\dim \Ker A^*$). The above characterization of badly approximable
functions can be easily generalized in the following way: if
$\f$ is a function in
$L^\be$ such that the essential norm
$\|H_\f\|_{\rm e}$ of the Hankel operator $H_\f$ (i.e., the {\it distance from
$H_\f$ to the set of compact operators}) is less than its norm, then $\f$
is badly approximable if and only if $\f$ has constant modulus almost  everywhere
on $\T$, $T_\f$ is Fredholm, and $\ind T_\f>0$ (see e.g., \cite{P}, Ch. 7, \S 5). 

Recall also that 
$$
\|H_\f\|=\dist_{L^\be}(\f,H^\be)\quad\mbox{and}\quad\|H_\f\|_{\rm e}=\dist_{L^\be}(\f,H^\be+C)
$$
(see, e.g., \cite{P}).

Let us proceed now to the case of matrix functions. We can consider the same problem
of finding a best analytic approximant for a given bounded function $\Phi$ with values
in the space $\mm_{m,n}$ of $m\times n$  matrices: for $\Phi\in L^\be(\mm_{m,n})$
find a bounded analytic $\mm_{m,n}$-valued function $F$ such that
$$
\|\Phi-F\|_{L^\be}=\dist_{L^\be(\mm_{m,n})}\Big(\Phi,H^\be(\mm_{m,n})\Big).
$$
Here 
$$
\|\Phi\|_{L^\be}\df\ess\sup_{\z\in\T}\|\Phi(\z)\|_{\mm_{m,n}},
$$
$\mm_{m,n}$ is equipped with the standard operator norm, and $H^\be(\mm_{m,n})$ is the space of
bounded analytic functions with values in $\mm_{m,n}$.

Again, it can be shown easily that a best approximant always exists. However, the situation
with uniqueness is quite different from the scalar case. Indeed,
suppose that $m=n=2$ and $u$ is a scalar continuous
badly approximable {\it unimodular function} (i.e., $|u(\z)|=1$ almost everywhere on $\T$).
Consider the matrix function
$\Phi=\left(\begin{array}{cc}u&\bs{0}\\\bs{0}&\0\end{array}\right)$. It is easy to see that for
any  scalar function $f$ in the unit ball of $H^\be$, the matrix function 
$\left(\begin{array}{cc}\0&\bs{0}\\\bs{0}&f\end{array}\right)$ is a best approximation of $\Phi$.

While it is possible to describe badly approximable matrix- and operator-valued functions, 
and we give such descriptions in this paper (the case of finite matrix functions was
treated in our earlier paper \cite{PT3}), this is not our main  goal. It turns out that in the
matrix case it is more natural to consider superoptimal approximations and \emph{very badly
approximable functions}. 

\medskip

{\bf Superoptimal approximations and very badly approximable matrix functions.}
Recall that for a matrix (or a bounded linear operator on Hilbert space)
$A$ the {\it singular values} $s_j(A)$, $j\ge0$, are defined by
$$
s_j(A)=\inf\{\|A-K\|:~\rank K\le j\}.
$$ 
Clearly, $s_0(A)=\|A\|$.

\medskip

{\bf Definition.}
Given a matrix function $\Phi\in L^\be(\mm_{m,n})$ we define inductively
the sets $\bs{\O}_j$, $0\le j\le\min\{m,n\}-1$, by
$$
\bs{\O}_0=\{F\in H^\be(\mm_{m,n})
:~F~\mbox{minimizes}~\ t_0\df\ess\sup_{\z\in\T}\|\Phi(\z)-F(\z)\|\};
$$
$$
\bs{\O}_j=\{F\in \O_{j-1}:~F~\mbox{minimizes}~\ 
t_j\df\ess\sup_{\z\in\T}s_j(\Phi(\z)-F(\z))\},\quad j>0.
$$
Functions in $\bigcap\limits_{k\ge 0} \bs{\O}_k = \bs{\O}_{\min\{m,n\}-1} $ are called {\it
superoptimal approximants} of $\Phi$ by bounded analytic matrix functions. The numbers
$t_j=t_j(\Phi)$ are called the {\it superoptimal singular values} of $\Phi$. 
Note that the functions
in $\bs{\O}_0$ are just the best approximants by bounded analytic matrix functions.

\medskip

As in the case of scalar functions, a bounded $m\times n$
matrix function $\Phi$ is called {\it badly approximable} if 
$$
\|\Phi\|_{L^\be}=\inf\{\|\Phi-F\|_{L^\be}:~F\in H^\be(\mm_{m,n})\}.
$$

We say that a matrix function $\Phi\in L^\be(\mm_{m,n})$ is called 
{\it very badly approximable} if the zero function $\0$ is a superoptimal approximant of $\Phi$.

The notion of superoptimal approximation can be extended to the case of
operator-valued functions. If $\h$ and $\K$ are Hilbert spaces, we denote
by $\B(\h,\K)$ the space of bounded linear operators from $\h$ to $\K$, 
$$
\B(\h)\df\B(\h,\h).
$$
We can identify an infinite-dimensional separable Hilbert
space with $\ell^2$ and identify operators on $\ell^2$ with infinite matrices.
Suppose that $\Phi\in L^\be(\B(\ell^2))$, i.e., $\Phi$ is a weakly measurable bounded
function that takes values in $\B(\ell^2)$, we can define the sequence
$\{\bs{\O}_j\}_{j\ge0}$ in the same way as for finite matrix functions. However, 
in the case of operator-valued functions we have to consider the infinite sequence
of the sets $\bs{\O}_j$. For $\Phi\in L^\be(\B(\ell^2))$, we say that a function
$F$ in $H^\be(\B(\ell^2))$ is a {\it superoptimal approximant of} $\Phi$ by bounded analytic
operator functions.

Badly approximable and very badly approximable infinite matrix functions can be defined in the
same way as in the case of finite matrix functions.

Note that if $\Phi$ is a matrix function of size $m\times\be$ or $\be\times n$, we can add to
$\Phi$ 
infinitely many zero rows or zero columns and reduce the problem to the case of matrix functions
of size
$\be\times\be$.

\medskip

{\bf The summary of earlier results.} First of all, let us mention 
that superoptimal approximation is more natural in the case of matrix or operator
functions because it is unique under mild natural assumptions on the function.
It was shown in \cite{PY} that if
$\Phi\in(H^\be+C)(\mm_{m,n})$ (i.e., all entries of
$\Phi$ belong to $H^\be+C$), then $\Phi$ has a unique superoptimal approximation $F$ by bounded
analytic matrix  functions. 
Moreover, it was shown in
\cite{PY} that
\bay
\label{sv}
s_j(\Phi(\z)-F(\z))=t_j(\Phi)\quad\mbox{for almost all}\quad\z\in\T.
\ey
Later this result was extended in \cite{T}, see also \cite{Pe1}, 
\cite{PT1} to operator-valued functions $\Phi$ for which the Hankel operator $H_\Phi$ is
compact.

The proof given in \cite{PY} was based on certain special factorizations (thematic
factorizations, see \S 4 of this paper for definitions). The approach in \cite{T} 
was more geometric and based on the notion of superoptimal weights.

The problem to describe the very badly approximable functions was posed in \cite{PY}. It follows from
\rf{sv} that if $\Phi$ is a very badly approximable function in \lb$(H^\be+C)(\mm_{m,n})$, then the singular values 
$s_j(\Phi(\z))$ are constant for almost all $\z\in\T$. Moreover, it was shown in \cite{PY} that if 
in addition to this $m\le n$ and
\mbox{$s_{m-1}(\Phi(\z))\ne0$} almost everywhere, 
then the Toeplitz operator \lb\mbox{$T_{z\Phi}:H^2(\C^n)\to H^2(\C^m)$}
has dense range (if $\Phi$ is a scalar function, the last condition is equivalent to the fact that 
$\ind T_\Phi>0$).
Note that the Toeplitz and the Hankel operators whose symbols are matrix functions can be 
defined in the same way as in the scalar case (see \rf{th}). Obviously, this necessary condition is equivalent to the 
condition $\Ker T_{\bar z\Phi^*}=\{\0\}$. In fact, the proof of necessity given in \cite{PY} allows one to obtain 
a more general result: if $\Phi$ is an arbitrary very badly approximable function in $(H^\be+C)(\mm_{m,n})$ and
$f\in\Ker T_{\bar z\Phi^*}$, then $\Phi^*f=\0$.

On the other hand, in \cite{PY} an example of a continuous $2\times2$ function $\Phi$ was given
such that
$s_0(\Phi(\z))=1$,
$s_1(\Phi(\z))=\a<1$, $\z\in\T$, $T_{z\Phi}$ is invertible but $\Phi$ is not even badly approximable.

The very badly approximable matrix functions of class $(H^\be+C)(\mm_{m,n})$ were characterized in
\cite{PY}  algebraically, in terms of so-called thematic factorizations. 

Later in \cite{PT} the above results of \cite{PY} were generalized to 
the broader context of matrix functions
$\Phi$ such that the essential norm $\|H_\Phi\|_{\rm e}$ of the Hankel 
operator $H_\Phi$ is less than the smallest
nonzero superoptimal singular value of $\Phi$. We call such matrix functions $\Phi$ {\it admissible}.
In particular, if $\Phi$ is an admissible very badly approximable $m\times n$ matrix function, then
the functions $s_j(\Phi(z))$ are constant almost everywhere on $\T$ and
$$
\Ker T_{\bar z\Phi^*}=\{f\in H^2(\C^n):~\Phi^*f=\0\}.
$$

In \cite{AP} another algebraic characterization of the set of very badly approximable 
admissible matrix functions was given in terms of canonical factorizations (see \S 5
for the definition).

We refer the reader to the book \cite{P}, which contains all the above information and 
results on superoptimal approximation and very badly approximable functions.

In \cite{PT3} we obtained a new criterion 
for an admissible matrix function to be very badly
approximable. In contrast with earlier criteria in terms of
certain special factorizations, it is more geometric and it is easier
to use it to verify whether a given matrix function is very badly approximable.
This criterion is given in terms of families of subspaces spanned by
Schmidt vectors of matrices $\Phi(\z)$, $\z\in\T$.

Recall that if $A$ is 
an $m\times n$ matrix and $s$ is a singular value of $A$, 
a nonzero vector $x\in\C^n$ is called
a {\it Schmidt vector} corresponding to $s$ if $A^*Ax=s^2x$.

Given a matrix function $\Phi$ in $L^\be(\mm_{m,n})$ and $\s>0$, we considered the subspace
${\frak S}_\Phi^{(\s)}(\z)$ of $\C^n$ spanned by the Schmidt vectors of
$\Phi(\z)$ that correspond to the singular values of $\Phi(\z)$ that are greater than or equal
to $\s$. The subspaces ${\frak S}_\Phi^{(\s)}(\z)$ are defined for almost all $\z\in\T$.
It was shown in \cite{PT3} that if $\Phi$ is an admissible very badly approximable
matrix functions, then for each $\s>0$, the family of subspaces ${\frak S}_\Phi^{(\s)}(\z)$,
$\z\in\T$, is {\it analytic}, i.e., there exist functions $g_1,\cdots,g_k$ in $H^2(\C^n)$
such that
\bay
\label{an}
{\frak S}_\Phi^{(\s)}(\z)=\spn\{g_1(\z),\cdots,g_k(\z)\}\quad\mbox{for almost all}\quad\z\in\T.
\ey
The same analyticity condition must also be imposed on the transposed function
$\Phi^{\rm t}$. However, it was shown in \cite{PT3} that the analyticity conditions on $\Phi$
and $\Phi^{\rm t}$ together with the earlier necessary conditions quoted above do not
guarantee that $\Phi$ is very badly approximable. 

However, it turned out that the above condition can be slightly modified to get a necessary and
sufficient condition. The main result of \cite{PT3} is the following theorem.

\medskip 

{\bf Theorem.} {\it Let $\Phi$ be an admissible matrix function. Then
$\Phi$ is  very badly approximable if and only if for each $\s>0$ equality {\em\rf{an}} holds
for functions $g_1,\cdots,g_k$ in $\Ker T_\Phi$. Moreover, this condition implies that $\Phi$
is very badly approximable even without the assumption that $\Phi$ is admissible.}

\medskip

Note that this condition in the case of a scalar function $\f$ means that $\f$ has constant
modulus and $\Ker T_\f\ne\{\0\}$, i.e., our criterion is a natural generalization of
the scalar results discussed above.

The uniqueness problem for superoptimal approximation of operator functions (infinite matrix
functions) was studied in \cite{T}, \cite{Pe1}, and \cite{PT1}. It was shown there that if 
the Hankel operator $H_\Phi$ is compact, then $\Phi$ has a unique superoptimal approximant by 
bounded analytic operator functions. In \cite{Pe1} and \cite{PT1} uniqueness was obtained with
the help of partial thematic factorizations (see \S 4 of this paper). We also refer the
reader to the monograph \cite{P} for the above results on superoptimal approximation of
operator functions.

\medskip

{\bf The purpose of this paper.} In this paper we study very badly approximable operator
functions. We consider the class of {\it admissible operator functions}. As in the case
of finite matrix functions, an operator function $\Phi$ is called admissible if
the essential norm $\|H_\Phi\|_{\rm e}$ of the Hankel operator $H_\Phi$ is less than each
nonzero superoptimal singular value of $\Phi$. 

In \S 4 we consider partial thematic factorizations of of admissible operator functions
(without the assumption of the compactness of $H_\Phi$ as it was done in
\cite{Pe1} and \cite{PT1}). In \S 5 we consider partial canonical factorizations
of operator functions.

The main result of the paper is 
a criterion of very bad approximability (Theorem 6.1) presented  in \S 6. 
It essentially says that
the theorem stated above also holds in the case of operator function. 
 
However, it turns out that
there is an important distinction between the case of finite matrix functions
and the case of infinite matrix functions. In the case of finite matrix functions
if $\Phi$ satisfies the hypotheses of the above theorem, then the zero function is the only
superoptimal approximant of $\Phi$. We show in this paper that in the case of infinite matrix
functions this is not true: {\it under the hypotheses of the above theorem $\Phi$ must be very
badly approximable, but it can have infinitely many superoptimal approximants}.

Note also that in the case of infinite matrix functions some proofs are considerably more
complicated than the proofs of the corresponding results for finite matrix functions
(e.g., the proofs of Theorems \ref{bal} and \ref{mv} given below).

In \S 2 we define inner, outer, and co-outer operator functions and prove
a theorem about inner-outer factorizations of co-outer operator functions.

In \S 3 we define balanced operator functions and prove that a inner and co-outer
function with finitely many columns has a balanced completion.

\

\section{\bf Inner and outer operator functions}
\setcounter{equation}{0}

\

In this section we define inner, outer, and co-outer operator functions and we prove that 
the inner factor in the inner-outer factorization of a co-outer function with finitely many
columns must also be co-outer.

Let $\h$ and $\K$ be separable Hilbert spaces. We denote by
$H^2_{\rm s}(\B(\h,\K))$ the space of analytic operator functions
$F$ that take values in the space of bounded linear operators form $\h$
to $\K$ and satisfy the following condition
$$
F(z)x\in H^2(\K)\quad\mbox{for every}\quad x\in\h.
$$
A function $F$ in $H^2_{\rm s}(\B(\h,\K))$ is called {\it inner} if $F(\z)$
is an isometric operator (i.e., $F(\z)^*F(\z)=I$) for almost all $\z\in\T$. A function 
$F$ in $H^2_{\rm s}(\B(\h,\K))$ is called {\it outer} if the set
$$
\{F q:~q~\mbox{ is a polynomial in }~H^2(\h)\}
$$
is dense in $H^2(\K)$.

It is well known (see e.g., \cite{N}) 
that each function $F$ in $H^2_{\rm s}(\B(\h,\K))$ admits an
inner-outer factorization, i.e., there exist an inner operator function $\Theta$
and an outer operator function $G$ such that $F=\Theta G$.

As we have mentioned in the introduction, we are going to identify operator functions with
infinite matrix functions. We say that an infinite matrix function $F$ 
is {\it co-outer} if the transposed function $F^{\rm t}$ is outer.

\begin{thm}
\label{cof}
Let $F$ be a co-outer operator function in $H^2_{\rm s}(\B(\C^d,\ell^2))$, $d<\be$. 
Suppose that
$$
F=\Theta G,
$$
where $\Theta$ is an inner operator function and $G$ is an outer
operator function. Then $\Theta$ is co-outer.
\end{thm}

\Pf Suppose that $\Theta\in H^\be(\B(\C^k,\ell^2))$ and $G\in H^2(\mm_{d,k})$.
Since $G$ is outer, it follows that $k\le d$. Suppose that 
$$
\Theta^{\rm t}=\co Q,
$$
where $\co$ is an inner matrix function and
$Q$ is an outer operator function. Since $\Theta$ is inner, it is easy to see that
$\co$ has size $k\times k$. We have
$$
\Theta=Q^{\rm t}\co^{\rm t}.
$$ 
Then
$$
F=Q^{\rm t}\co^{\rm t}G,
$$
and so by the hypotheses of the theorem,
$$
F^{\rm t}=G^{\rm t}\co Q
$$
is an outer function. It follows that $G^{\rm t}$ must be outer, and so $k=d$.
Clearly, $G^{\rm t}\co H^\be(\mm_{d,d})$ must be dense in $H^2(\mm_{d,d})$. However, the
determinants of all matrix functions in $G^{\rm t}\co H^\be(\mm_{d,d})$ must be divisible
by $\det\co$ which is a scalar inner function. Thus $\det\co$ is constant, and so
$\co^*=\co^{-1}\in H^\be(\mm_{d,d})$ which implies that $\co$ is constant, and so 
$\Theta^{\rm t}$ is co-outer. $\bl$

\

\section{\bf Balanced matrix functions}
\setcounter{equation}{0}

\

In this section we introduce the notion of balanced unitary-valued functions
and prove the existence of balanced completions for inner and co-outer functions that have
finitely many columns.

\medskip

{\bf Definition.} A {\it balanced infinite matrix function} is a 
unitary-valued matrix function 
of the form $\left(\begin{array}{cc}\U&\ov\Theta\end{array}\right)$, where
$\U$ and $\Theta$ inner and co-outer matrix
functions. 

If $\U$ has $r$ columns, we say that the function
$\left(\begin{array}{cc}\U&\ov\Theta\end{array}\right)$ is {\it $r$-balanced}.
1-balanced functions are also called {\it thematic matrix functions}.

\medskip

We are going to prove that an inner matrix function with finitely many columns can 
be completed to a balanced matrix function.

\medskip

Let $r$ be a positive integer and let
$\U$ be an inner matrix function in $H^\be(\C^r,\ell^2)$.
Consider the subspace $\cL\df\Ker T_{\U^t}$ of $H^2(\ell^2)$.
Clearly, it is invariant under multiplication by $z$, and
so there exists an inner matrix function $\Theta$ such that 
$\cL=\Theta H^2(\K)$, where $\K=\ell^2$ or $\K=\C^m$ for some positive $m$.
The proof of the following theorem in the special case $r=1$ can be found in \cite{P}, 
Ch. 14, \S 18. In the general case the proof is algebraically more complicated. 
Note that a close result was obtained in \cite{C}, see also \cite{H}, Lect. IX.

\begin{thm}
\label{bal}
Let $\U$ and $\Theta$ be as above.
Then $\Theta$ is co-outer and the matrix function
$\left(\begin{array}{cc}\U&\ov\Theta\end{array}\right)$ is unitary-valued.
\end{thm}

Before proceeding to the proof, we introduce a notion. Let 
$$
A=\left(\begin{array}{cccc}
a_{11}&a_{12}&\cdots&a_{1r}\\
a_{21}&a_{22}&\cdots&a_{2r}\\
\vdots&\vdots&\ddots&\vdots\\
a_{r+1~1}&a_{r+1~2}&\cdots&a_{r+1~r}
\end{array}\right)
$$
be an $(r+1)\times r$ matrix. For $1\le j\le r+1$, we put
$$ 
\a_j=(-1)^{j}\det\left(\begin{array}{ccc}
a_{11}&\cdots&a_{1r}\\
\vdots&\ddots&\vdots\\
a_{j-1~1}&\cdots&a_{j-1~r}\\
a_{j-1~1}&\cdots&a_{j-1~r}\\
\vdots&\ddots&\vdots\\
a_{r+1~1}&\cdots&a_{r+1~r}
\end{array}\right).
$$
In other words, we multiply $(-1)^j$ by the minor obtained from $A$ by deleting the
$j$th row. The vector $A_{\rm ass}\df\{\a_j\}_{1\le j\le r+1}$ is called the {\it vector associated with}
$A$. 

\medskip

\Pf The proof of the fact that $\Theta$ is co-outer is exactly the same
as in the case $r=1$, see \cite{P}, Ch. 14, Lemma 18.3. Let us show that 
$\left(\begin{array}{cc}\U&\Theta\end{array}\right)$ is unitary-valued.
The fact that $\left(\begin{array}{cc}\U&\Theta\end{array}\right)$ takes isometric values
almost everywhere on $\T$ follows immediately from the definition of $\Theta$.
To prove that it is unitary-valued, it suffices to show that 
$\dim\Ker\Theta^{\rm t}(\z)\le r$ for almost all $\z\in\T$.

Let
$$
\U=\left(\begin{array}{cccc}\up_{01}&\up_{02}&\cdots&\up_{0r}\\
\up_{11}&\up_{12}&\cdots&\up_{1r}\\
\up_{21}&\up_{22}&\cdots&\up_{2r}\\
\vdots&\vdots&\ddots&\vdots
\end{array}\right).
$$

Clearly, the matrix function $\U$ has rank $r$ almost everywhere on $\T$. Without loss of generality
we may assume that
\bay
\label{det}
\det\left(\begin{array}{cccc}\up_{01}&\up_{02}&\cdots&\up_{0r}\\
\up_{11}&\up_{12}&\cdots&\up_{1r}\\
\up_{21}&\up_{22}&\cdots&\up_{2r}\\
\vdots&\vdots&\ddots&\vdots\\
\up_{r-1~1}&\up_{r-1~2}&\cdots&\up_{r-1~r}
\end{array}\right)\ne\0.
\ey

Consider the bounded analytic matrix function $G$ defined in the following way:
$$
G=\left(\begin{array}{cccccc}
\a_0^{[0]}&\a_0^{[1]}&\a_0^{[2]}&\a_0^{[3]}&\cdots\\[.2cm]
\a_1^{[0]}&\a_1^{[1]}&\a_1^{[2]}&\a_1^{[3]}&\cdots\\[.2cm]
\vdots&\vdots&\vdots&\vdots&\ddots\\[.2cm]
\a_{r-1}^{[0]}&\a_{r-1}^{[1]}&\a_{r-1}^{[2]}&\a_{r-1}^{[3]}&\cdots\\[.2cm]
\a_r^{[0]}&\0&\0&\0&\cdots\\[.2cm]
\0&\a_r^{[0]}&\0&\0&\cdots\\[.2cm]
\0&\0&\a_r^{[0]}&\0&\cdots\\[.2cm]
\vdots&\vdots&\ddots&\ddots&\ddots
\end{array}\right),
$$
where for $k\ge0$, the $H^\be$ functions $\a^{[k]}_m$, $0\le m\le r$, are the components of the vector
function $A^{[k]}_{\rm ass}$ associated with the matrix function $A^{[k]}$ defined by
$$
A^{[k]}=\left(\begin{array}{cccc}
\up_{01}&\up_{02}&\cdots&\up_{0r}\\
\up_{11}&\up_{12}&\cdots&\up_{1r}\\
\up_{21}&\up_{22}&\cdots&\up_{2r}\\
\vdots&\vdots&\ddots&\vdots\\
\up_{r-1~1}&\up_{r-1~2}&\cdots&\up_{r-1~r}\\
\up_{r+k~1}&\up_{r+k~2}&\cdots&\up_{r+k~r}
\end{array}\right).
$$
Note that $\a_r^{[0]}$ is nothing but the determinant on the left-hand side of \rf{det}.

It is an elementary exercise in linear algebra to verify that $\U^{\rm t}G=\0$. It follows that
$G$ admits a factorization $G=\Theta Q$, where $Q$ is an $H^\be$ matrix function. Hence, 
to verify that $\dim\Ker\Theta^{\rm t}(\z)\le r$, it suffices to show that
$\dim\Ker G^{\rm t}(\z)\le r$. Recall that by \rf{det}, $\a_r^{[0]}(\z)\ne0$ 
for almost all $\z\in\T$. 
Assume that $\z\in\T$ and $\a_r^{[0]}(\z)\ne0$.
Suppose that the vector $x=\{x_j\}_{j\ge0}$ belongs to $\Ker G^{\rm t}(\z)$.
If we look at the $r$th coordinate of the vector $G^{\rm t}(\z)x$, 
we observe that $\a_r^{[0]}(\z)x_r$
is uniquely determined by $x_0,x_1,\cdots,x_{r-1}$. 
Since $\a_r^{[0]}(\z)\ne0$, it follows that
$x_r$ is uniquely determined by $x_0,x_1,\cdots,x_{r-1}$. 
If we look now at the next component of the vector
$G^{\rm t}(\z)x$, we observe that $x_{r+1}$  
is uniquely determined by $x_0,x_1,\cdots,x_r$, etc.
This completes the proof. $\bl$

\

\section{\bf Partial thematic factorizations}
\setcounter{equation}{0}

\

In the case when $\Phi$ is an operator function such that the Hankel operator $H_\Phi$ is
compact and $F\in\bs{\O}_d$, partial thematic factorizations of $\Phi-F$ were
constructed in \cite{Pe1}. In this section we consider the more general case when
$\Phi$ is an admissible operator function. 

Suppose that $\Phi$ is function in $L^\be(\B(\ell^2))$ such that
\bay
\label{assum}
\|H_\Phi\|_{\rm e}<\|H_\Phi\|
\ey
and $F\in H^\be(\B(\ell^2))$ is a best approximant of $\Phi$.
Then $H_\Phi$ has a maximizing vector $f$, the function $g=\|H_\Phi\|^{-1}\bar z\ov{H_\Phi f}$
is a maximizing vector of $H_{\Phi^{\rm t}}$. The functions $f$ and $g$ admit factorizations
$$
f=\vt_1hv,\quad g=\vt_2hw,
$$
where $h$ is a scalar outer function, $\vt_1$ and $\vt_2$ are scalar inner functions, and $v$
and $w$ are inner and co-outer column functions. 

By Theorem \ref{bal}, the column functions $v$ and $w$ have thematic (1-balanced) completions:
$$
V=\left(\begin{array}{cc}v&\ov{\Theta}\end{array}\right)\quad\mbox{and}\quad
W^{\rm t}=\left(\begin{array}{cc}w&\ov{\Xi}\end{array}\right)
$$
The function $\Phi-F$ admits the following factorization:
$$
\Phi-F=W^*\left(\begin{array}{cc}t_0u&\0\\\0&\Psi\end{array}\right)V^*,
$$
where $u=\bar z\bar\vt_1\bar\vt_2\bar h/h$ and $\|\Psi\|_{L^\be}\le t_0=t_0(\Phi)=\|H_\Phi\|$
(see \cite{P}, Ch. 14, \S 18). Moreover, under the assumption \rf{assum}, $T_u$ is Fredholm and
$\ind T_u>0$.

Such factorizations are called {\it partial thematic factorizations of order} 1. 

As in the case of finite matrix functions (see \cite{PT} or \cite{P}, Ch. 14, \S 4) the
following crucial inequality holds:
\bay
\label{essnorm}
\|H_\Psi\|_{\rm e}\le\|H_\Psi\|_{\rm e}.
\ey
Another important result that can be established in the same way as in the
case of finite matrix functions is that under the assumption \rf{assum}
the operator functions $\Theta$ and $\xi$ are left-invertible in $H^\be$ 
(see \cite{PT} or \cite{P}, Ch. 14, \S 4).

A function $\Phi$ satisfying \rf{assum} is badly approximable if and only if
it admits a partial thematic factorization of order 1:
$$
\Phi=W^*\left(\begin{array}{cc}t_0u&\0\\\0&\Psi\end{array}\right)V^*
$$ 
(the part ``if'' holds even without the assumption \rf{assum}). Moreover, $\Phi$ is 
very badly approximable if and only if $\Psi$ is very badly approximable.

If $\Phi$ is admissible and $H_\Psi\ne0$, due to inequality \rf{essnorm} 
we can apply the same procedure to $\Psi$. If $F\in\bs{\O}_1$, then
$\Phi-F$ admits a thematic factorization of order 2, i.e.,
$$
\Phi-F=W^*\left(\begin{array}{cc}1&\0\\\0&W_1^*\end{array}\right)
\left(\begin{array}{ccc}t_0u_0&\0&\0\\\0&t_1u_1&\0\\
\0&\0&\L\end{array}\right)
\left(\begin{array}{cc}1&\0\\\0&V_1^*\end{array}\right)V^*,
$$
where $V,\,V_1,\,W^{\rm t},\,W_1^{\rm t}$ are thematic operator functions,
$u_0$ and $u-1$ are scalar very badly approximable functions such that
$\|H_{u_j}\|_{\rm e}<1$, and $\|\Psi\|_{L^\be}\le t_1$.

If $\Phi$ is admissible, we can continue this process and obtain partial thematic factorization
of an arbitrary order. 

In particular an admissible operator function $\Phi$ is very badly approximable if and only 
is for each positive integer $r$ it admits a partial thematic factorizations of order $r$.

\

\section{\bf Partial canonical factorizations}
\setcounter{equation}{0}

\

As in the case of finite matrix functions (see \cite{PT3}), to obtain to obtain a geometric
characterization of very badly approximable operator functions, it is more important to deal
with canonical factorizations rather than with thematic factorizations.

\begin{thm}
\label{mv}
Let $\Phi$ be a matrix function in $L^\be(\B(\ell^2))$ such that
that \lb$\|H_\Phi\|_{\rm e}<\|H_\Phi\|$ and let
$r$ be the multiplicity of the superoptimal singular value $t_0(\Phi)$.
Suppose that $\M$ is the minimal shift invariant subspace
of $H^2(\ell^2)$ that contains all maximizing vectors of $H_\Phi$. Then
$$
\M=\U H^2(\C^d),
$$
where $\U$ is an inner and co-outer function of size $\be\times r$.
\end{thm}

\Pf Since $\M$ is shift invariant, it has the form
$$
\M=\U H^2(\K),
$$
where $\K$ is a separable Hilbert space and $\U$ is an inner operator function. 
Since $\|H_\Phi\|_{\rm e}<\|H_\Phi\|$, it is easy to see that the space of maximizing
vectors of $H_\Phi$ is finite-dimensional, and so $\dim\K<\be$. Put $d=\dim\K$ and 
$\K=\C^d$.

Let us show that $d\ge r$. In \cite{PY} (see also Lemma 1.2 of \cite{PY2}) in the case of
finite matrix functions of class $H^\be+C$ a finite sequence 
$$
f^{(0)}_1,\cdots,f^{(0)}_{k_0},f^{(1)}_1,\cdots,f^{(1)}_{k_1},\cdots,
f^{(r-1)}_1,\cdots,f^{(r-1)}_{k_{r-1}}
$$
of maximizing vectors of $H_\Phi$ was constructed. It is easy to verify that it has the
following property:
\bay
\label{prop}
\max_{\z\in\dd}\,\dim\,\spn\left\{f^{(j)}_k(\z):~0\le j\le r-1,~1\le k\le k_j\right\}=r.
\ey
This construction was generalized in \cite{PT1} to the case of
finite matrix functions $\Phi$ satisfying the condition $\|H_\Phi\|_{\rm e}<\|H_\Phi\|$
and in \cite{PT} to the case of infinite matrix functions $\Phi$ such that $H_\Phi$
is compact (see also Chap. 14 of \cite{P}). It can easily be verified that exactly the
same construction also works in the case of infinite matrix functions $\Phi$ satisfying the
condition $\|H_\Phi\|_{\rm e}<\|H_\Phi\|$ and \rf{prop} holds. It follows immediately from
\rf{prop} that $d\ge r$.

Let us now show that $d\le r$. 
Let $F$ be a function in $\bs{\O}_r$. Consider a partial canonical factorization 
of $\Phi-F$. It has the form
$$
\Phi-F={\frak W}\left(\begin{array}{cccccc}
t_0{\frak U}&\0\\
\0&\Psi
\end{array}\right){\frak V},
$$
where ${\frak W}$ and ${\frak V}$ are infinite unitary-valued functions, ${\frak U}$ is
an $r\times r$ unitary-valued function, and $\|\Psi\|_{L^\be}=t_r<t_0$. It follows that
the subspace spanned by the maximizing vectors of $(\Phi-F)(\z)$ has dimension $r$
for almost all $\z\in\T$. 

For every function $f\in\M$ the vector $f(\z)$ is a maximizing
vector of $\Phi(\z)$ for almost all $\z\in\T$ (see Lemma 15.2 in Ch. 14 of 
\cite{P}, note that in \cite{P} the result is stated for finite matrix functions, but the proof
given there works for infinite matrix functions too). It is easy to see now that if $d>r$, 
then the subspace spanned by the maximizing vectors of $(\Phi-F)(\z)$ has dimension 
at least $d$. 

It remains to show that $\U$ is co-outer. Without loss of generality we may assume that
$\|\Phi\|_\be=\|H_\Phi\|=1$.
Consider the subspace of $H^2(\ell^2)$ spanned by the maximizing vectors of $H_\Phi$.
It must be finite-dimensional. Let $f_1,\cdots,f_s$ be a basis of this subspace and
let $F$ be the matrix function whose columns are $f_1,\cdots,f_s$. Consider the
inner--outer factorization of $F^{\rm t}$:
$$
F^{\rm t}=\co G,
$$
where $\co$ is an inner matrix function of size $s\times k$, $k\ge s$,
and $G$ is an outer matrix function
of size $k\times\be$. Then
$$
F=G^{\rm t}\co^{\rm t},
$$
and so 
\bay
\label{O}
G^{\rm t}=G^{\rm t}\co^{\rm t}\ov{\co}=F\ov{\co}.
\ey

Since the functions $f_j$ are maximizing vectors of $H_\Phi$, it follows that
for almost all $\z\in\T$, the vectors $f_j(\z)$ are maximizing vectors of $\Phi(\z)$ and 
$\Phi f_j\in H^2_-(\ell^2)$ (see \cite{P}, Theorem 2.3 of Ch. 2). Thus for almost all $\z\in\T$,
the restriction of $\Phi(\z)$ to $\Range F(\z)$ is an isometry and 
$\Phi F\in H^2_-\big(\B(\C^s,\ell^2)\big)$. It follows that
$$
\Phi G^{\rm t}=\Phi F\ov{\co}\in H^2_-\big(\B(\C^k,\ell^2)\big).
$$
Suppose now that $g$ is a column of $G^{\rm t}$. Then $\Phi g\in H^2_-(\ell^2)$. Since 
$$
g(\z)\in\Range G^{\rm t}(\z)\subset\Range F(\z)\quad\mbox{for almost all }~\z\in\T,
$$ 
it follows that
$$
\|\Phi(\z)g(\z)\|_{\ell^2}=\|g(\z)\|_{\ell^2}\quad\mbox{almost everywhere on }~\T.
$$
Thus
$$
\|H_\Phi g\|=\|\pp_-\Phi g\|=\|\Phi g\|=\|g\|,
$$
and so all columns of $G^{\rm t}$ are maximizing vectors of $H_\Phi$. Since the columns
of $F$ form a basis in the space of maximizing vectors, it follows from \rf{O}
that $\co$ is a constant isometric matrix.

Clearly, $\M$ is the minimal invariant subspace of multiplication by $z$ on $H^2(\ell^2)$
that contains the columns of $F$. Consider the subspace minimal invariant subspace
$\M_1$ that contains the columns of $G^{\rm t}$. Since $F=G^{\rm t}\co^{\rm t}$ and 
$\co$ is a constant matrix, it follows that $\M\subset\M_1$. On the other hand, the columns of
$G^{\rm t}$ are maximizing vectors of $H_\Phi$, and so $\M_1\subset\M$. Thus $\M_1=\M$.

Now it is easy to see that $\U$ is inner factor of the inner-outer factorization of $G^{\rm t}$.
It follows now from Theorem \ref{cof} that $\U$ is co-outer. $\bl$

Consider now the matrix function $\Phi^{\rm t}$. Let $\N$ be the shift-invariant subspace of
$H^2(\ell^2)$ spanned by the maximizing vectors of $H_{\Phi^{\rm t}}$.
Then by Theorem \ref{mv}, $\N$ has the form $\O H^2(\C^r)$, where $\O$ is an inner and
co-outer matrix function. By Theorem \ref{bal}, there exist inner and co-outer
matrix functions $\Theta$ and $\Xi$ such that
\bay
\label{VW}
\V\df\left(\begin{array}{cc}\U&\ov{\Theta}\end{array}\right)\quad\mbox{and}\quad
\W^{\text t}\df\left(\begin{array}{cc}\O&\ov{\Xi}\end{array}\right)
\ey
are unitary-valued matrix functions.

The proof of the following result is exactly the same as the proof of Theorem
15.3 of Ch. 14 of \cite{P} for finite matrix functions (see also \cite{AP}).

\begin{thm}
\label{fact}
Let $\Phi$ be a function in $L^\be(\B(\ell^2))$ such that 
{\em$\|H_\Phi\|_{\text e}<t_0=\|H_\Phi\|$}.  Let $r$ be the number of superoptimal singular
values of $\Phi$ equal to 
$t_0$. Suppose that $F$ is a best approximation
of $\Phi$ by analytic matrix functions. Then $\Phi-F$ admits a factorization
of the form
\bay
\label{chast}
\Phi-F=\W^*\left(\begin{array}{cc}t_0U&\0\\\0&\Psi
\end{array}\right)\V^*,
\ey
where $\V$ and $\W$ are given by {\em\rf{VW}}, $U$ is an $r\times r$ 
unitary-valued very badly \lb approximable matrix function such that 
{\em$\|H_U\|_{\text e}<1$}, and $\Psi$ is a \lb matrix-function in
$L^\be(\B(\ell^2))$ such that \mbox{$\|\Psi\|_{L^\be}\le t_0$} and
\lb$\|H_\Psi\|=t_r(\Phi)<\|H_\Phi\|$. Moreover, $U$ is uniquely determined by the 
choice of $\U$ and $\O$ and does not depend on the choice of $F$.
\end{thm}

As in the case of finite matrix functions, under the hypotheses of Theorem \ref{fact}
the following inequality holds
\bay
\label{esnorm}
\|H_\Psi\|_{\rm e}\le\|H_\Phi\|_{\rm e}.
\ey
it can be deduced from \rf{essnorm} in exactly the same way as in \cite{AP}
(see also Theorem 15.12 of Ch. 15 of \cite{P}).

Moreover, under the hypotheses of Theorem \ref{fact}, the operator functions $\Theta$ and $\Xi$
in \rf{VW} are left-invertible in $H^\be$. Again, this can be deduced from the same
results for partial thematic factorizations (see \S 4) in the same way it was done in the case
of finite matrix functions in \cite{AP} (see also \cite{P}, Ch. 14, \S 5). This
left-invertibility property of $\Theta$ and $\Xi$ is important in the main result
of the next section.

The following theorem can be considered as a converse of Theorem \ref{fact}. These
two theorems together give a characterization of the badly approximable matrix
functions $\Phi$ satisfying the condition $\|H_\Phi\|_{\text e}<t_0=\|H_\Phi\|$. 
Note however that
we do not need this condition to prove that functions that admit a factorization of the form
\rf{chast}. Moreover, in the following theorem we can also relax the assumptions on $U$ imposed
in Theorem \ref{fact}.

\begin{thm}
\label{plokho}
Let $\Phi$ be an infinite matrix function of the form
\bay
\label{ordre1}
\Phi=\W^*\left(\begin{array}{cc}\s U&\0\\\0&\Psi
\end{array}\right)\V^*,
\ey
where $\s>0$, $\V$ and {\em$\W^{\text t}$} are $r$-balanced matrix functions, 
$U$ is an $r\times r$ unitary-valued matrix function such that
the shift-invariant subspace of $H^2(\C^r)$ spanned by the
maximizing vectors of $H_U$ coincides with $H^2(\C^r)$,
and $\|\Psi\|_\be\le\s$. Then
$\Phi$ is badly approximable and $t_0(\Phi)=\cdots=t_{r-1}(\Phi)=\s$.
Moreover, $\Phi$ is very badly approximable if and only if
$\Psi$ is very badly approximable.
\end{thm}

The proof of Theorem \ref{plokho} is exactly the same as the proof of
Theorem 15.7 of Ch. 14 of \cite{P} for finite matrix functions (see also \cite{AP}).

Consider now the sequence
$$
t_0=\cdots=t_{r_1-1}>t_{r_1}=\cdots=t_{r_2-1}>\cdots>
t_{r_{\iota-1}}=\cdots=t_{r_\iota-1}>\cdots
$$
of superoptimal singular values of $\Phi$. Let
$$
\s_0>\s_1>\s_2>\cdots
$$
be the sequence of distinct superoptimal singular values of $\Phi$, i.e.,
$$
\s_0=t_0=\cdots=t_{r_1-1},\quad\s_1=t_{r_1}=\cdots=t_{r_2-1},\quad\mbox{etc}.
$$

If $\|H_\Phi\|_{\rm e}<\s_1$, we can apply Theorem \ref{fact} to the matrix function $\Psi$.
Now if $\|H_\Phi\|_{\rm e}<\s_2$, then by \rf{esnorm}, $\|H_\Psi\|_{\rm e}<\s_2$,
and so we can continue this process and obtain the following
result in
exactly the same way as in the case of finite matrix functions in \cite{P}, Ch. 14, \S 15.

\begin{thm}
\label{pcf}
Let $\Phi$ be a function in $L^\be(\B(\ell^2))$ such that 
$\|H_\Phi\|_{\rm e}<\s_{d-1}$. Let $F$ be an arbitrary matrix function in 
$\bs{\O}_{r_d}$. Then $\Phi-F$ admits a factorization
\begin{align}
\label{partial}
\Phi-F
=&\W_0^*\cdots\W^*_{d-1}
\left(\begin{array}{ccccc}\s_0U_0&\0&\cdots&\0&\0\\\0&\s_1U_1&\cdots&\0&\0\\
\vdots&\vdots&\ddots&\vdots&\vdots\\
\0&\0&\cdots&\s_{d-1}U_{d-1}&\0\\
\0&\0&\cdots&\0&\Psi\end{array}\right)
\V^*_{d-1}\cdots\V_0^*,
\end{align}
where the $U_j$ are $(r_{j+1}-r_j)\times(r_{j+1}-r_j)$ very badly approximable 
unitary-valued functions such that $\|H_{U_j}\|_{\text e}<1$,
$$
\V_j=\left(\begin{array}{cc}\bs{I}_{r_j}&\0\\\0&\breve{\V}_j\end{array}\right)
\quad\mbox{and}\quad
\W_j=\left(\begin{array}{cc}\bs{I}_{r_j}&\0\\\0&\breve{\W}_j\end{array}\right),
\quad1\le j\le d-1,
$$
$\breve{\V}_j$ and $\breve{\W}^{\rm t}_j$ are $(r_{j+1}-r_j)$-balanced matrix 
functions, and $\Psi$ is a matrix function satisfying
$$
\|\Psi\|_{L^\be}\le t_{r_{d-1}},\quad\mbox{and}\quad 
\|H_{\Psi}\|<t_{r_{d-1}}.
$$
\end{thm}

Factorizations of the form \rf{partial}
with the $\s_j$, $U_j$, $\V_j$, and $\W_j$ as in Theorem \ref{pcf}
are called {\it partial canonical factorizations 
(or partial canonical factorizations of order $d$)}.

Now we can state the following description of very badly approximable
matrix functions.

\begin{thm}
\label{vbaf}
Let $\Phi$ be an admissible function in $L^\be(\B(\ell^2))$.
If $\Phi$ is very badly approximable, then for each $d$ with nonzero $\s_{d-1}$
the matrix function $\Phi$ admits a partial canonical factorization of the form
{\em\rf{fact}}. 
\end{thm}

The proof of Theorem \ref{vbaf} is exactly the same as in the case of finite matrix functions,
see \cite{P}, Ch. 14, \S 15 (see also \cite{AP}). Finally, we state the converse of
Theorem \ref{vbaf}, which is valid without the admissibility assumption.

\begin{thm}
\label{converse}
Let $\Phi$ be a function in $L^\be(\B(\ell^2))$ such that $\Phi$ admits a
partial canonical factorization of the form {\em\rf{partial}} whenever $\s_{d-1}>0$.
Then $\Phi$ is very badly approximable and
$$
t_\varkappa(\Phi)=\left\{\begin{array}{ll}\s_0,&\varkappa<r_1,\\
\s_j,&r_j\le \varkappa<r_{j+1}.
\end{array}\right.
$$
\end{thm}

\medskip

{\bf Example.} As we have mentioned in the Introduction there is an important difference
between the case of finite matrix functions and the case of infinite matrix functions. In the
case of finite matrix functions the hypotheses of Theorem \ref{vbaf} guarantee that the zero
function the only superoptimal approximant. It turns out that in the case of infinite matrix
functions this is not true. Consider the following example.

Let $\{u_j\}_{j\ge0}$, be a sequence of scalar badly approximable functions such that
$$
|u_j(\z)|=1\quad\mbox{for almost all}\quad\z\in\T\quad\mbox{and}\quad\|H_{u_j}\|_{\rm e}<1
$$
and let $\{t_j\}_{j\ge0}$ be a decreasing sequence of positive numbers such that
$$
\lim_{j\to\be}t_j>0.
$$
Consider the infinite matrix function 
$$
\Phi=\left(\begin{array}{ccccc}\0&\0&\0&\0&\cdots\\
\0&t_0u_0&\0&\0&\cdots\\
\0&\0&t_1u_1&\0&\cdots\\
\0&\0&\0&t_2u_2&\cdots\\
\vdots&\vdots&\vdots&\vdots&\ddots\\
\end{array}\right).
$$
Obviously, for every $d\in\Z_+$, there is constant unitary matrix $V_d$ such that
\bay
\label{post}
\Phi=V_d^*
\left(\begin{array}{cccccccc}
t_0u_0&\0&\0&\cdots&\0&\0&\cdots\\
\0&t_1u_1&\0&\cdots&\0&\0&\cdots\\
\vdots&\vdots&\ddots&\vdots&\vdots&\vdots&\vdots\\
\0&\0&\cdots&t_{d-1}u_{d-1}&\0&\0&\cdots\\
\0&\0&\cdots&\0&\0&\0&\cdots\\
\0&\0&\cdots&\0&\0&t_du_d&\cdots\\
\vdots&\vdots&\vdots&\vdots&\vdots&\vdots&\ddots
\end{array}\right)
V_d.
\ey
Clearly, the right-hand side of \rf{post} is a partial canonical factorization of $\Phi$, and so
by Theorem \ref{converse}, $\Phi$ is very badly approximable and $t_j(\Phi)=t_j$

On the other hand, if $f$ is an arbitrary scalar function in $H^\be$ with 
$$
\|f\|_\be\le\lim_{j\to\be}t_j
$$
and 
$$
F=\left(\begin{array}{ccccc}f&\0&\0&\cdots\\
\0&\0&\0&\cdots\\
\0&\0&\0&\cdots\\
\vdots&\vdots&\vdots&\ddots\\
\end{array}\right),
$$
then 
$$
\left(\begin{array}{ccccc}-f&\0&\0&\0&\cdots\\
\0&t_0u_0&\0&\0&\cdots\\
\0&\0&t_1u_1&\0&\cdots\\
\0&\0&\0&t_2u_2&\cdots\\
\vdots&\vdots&\vdots&\vdots&\ddots\\
\end{array}\right)
$$
and since, obviously,
$$
s_j\big((\Phi-F)(\z)\big)=t_j,\quad j\in\Z_+,~\z\in\T,
$$
it follows that $F$ is a superoptimal approximant of $\Phi$.

\medskip

To make the conclusion that an admissible infinite matrix function
has a unique superoptimal approximant, we need the condition that
\bay
\label{nol'}
\lim_{j\to\be}t_j(\Phi)=0.
\ey
Indeed, if $F_1$ and $F_2$ belong to $\bs{\O}_r$, then we can consider
partial thematic factorizations of $\Phi-F_1$ and $\Phi-F_2$ and
see that
$$
\|F_1-F_2\|_\be=\big\|(\Phi-F_1)-(\Phi-F_2)\big\|_\be\le2t_r(\Phi).
$$
In particular, if both $F_1$ and $F_2$ are superoptimal approximants, then by 
\rf{nol'}, $F_1=F_2$. However, if $\Phi$ is admissible and satisfies \rf{nol'},
then $\lim\limits_{j\to\be}s_j(T_\Phi)=0$, and so $H_\Phi$ is compact.

\

\section{\bf Very badly approximable functions}
\setcounter{equation}{0}

\

In this sections we obtain a necessary and sufficient condition for
an admissible infinite matrix function to be very badly approximable.
Let $\Phi\in L^\be(\B(\ell^2))$. Put
$$
t_\be(\Phi)=\lim_{j\to\be}t_j(\Phi).
$$

As in the case of finite matrix functions, for $\s>t_\be(\Phi)$, we consider the subspace
${\frak S}_\Phi^{(\s)}(\z)$ that is the linear span of the Schmidt vectors
of $\Phi(\z)$ that correspond to the singular values of $\Phi(\z)$ that 
are greater than or equal to $\s$.
The subspaces ${\frak S}_\Phi^{(\s)}(\z)$ are defined for almost all $\z\in\T$.

\medskip

{\bf Definition.} Let $L(\z)$, $\z\in\T$, be a family of subspaces of 
$\ell^2$  that is defined almost everywhere on $\T$. 
We say that functions $\xi_1,\cdots,\xi_l$ in $H^2(\ell^2)$ {\it span the family} $L$ if
$L(\z)=\spn\{\xi_j(\z):~1\le j\le l\}$ for almost all $\z\in\T$.

\medskip 

We consider in this section the following condition:

\begin{enumerate}
\item[(C)] {\em for each $\s>t_\be(\Phi)$, the family of subspaces 
${\frak S}_\Phi^{(\s)}$ is analytic and spanned by finitely many functions in
$\Ker T_\Phi$}.
\end{enumerate}

As in the case of finite matrix functions (see \cite{PT3}), it is easy to see that condition (C)
implies that the functions $\z\mapsto s_j\big(\Phi(\z)\big)$, $j\in\Z_+$, are constant almost
everywhere on $\T$.

The following theorem is the main result of this section.

\begin{thm}
\label{plo}
If $\Phi$ is an admissible very badly approximable matrix function in $L^\be(\B(\ell^2))$,
then $\Phi$ satisfies {\em(C)}.

Conversely, if $\Phi$ is an arbitrary function in 
$L^\be(\B(\ell^2))$ that satisfies  {\em(C)}, then $\Phi$ is very badly approximable.
\end{thm}

{\bf Remark.} As we have already mentioned, there is an important difference between 
the case of finite matrix functions and the case of infinite matrix functions.
In the case of finite matrix functions condition (C) also implies that the zero function is the
only superoptimal approximant. In the case of infinite matrix functions this is not true. 
Indeed, it is easy to see that the matrix function given in the example at the end of the
previous section satisfies condition (C). However, it has infinitely many 
superoptimal approximants.

\medskip

The necessity of condition (C) can be obtained from Theorem \ref{vbaf} in exactly the same
way as it was done in \cite{PT3}, Theorem 4.1 in the case of finite matrix functions. On the
other hand, the proof of the sufficiency of (C) given in \cite{PT3} works only for finite
matrices. It has to be slightly modified to work
in the case of infinite matrix functions. 

Here we present a proof based on canonical factorization. 
Note that the proof based on superoptimal weights that was presented in \S 5 of \cite{PT3} in
the case of finite matrix functions also works (with obvious modifications). 

\medskip

{\bf Proof of the sufficiency of (C).} Suppose that $\Phi$ satisfies (C). As we have already
observed, the functions $\z\mapsto s_j\big(\Phi(\z)\big)$, $j\in\Z_+$, are constant almost
everywhere on $\T$. Let
\bay
\label{sigma}
\s_0>\s_1>\s_2>\cdots
\ey
be positive numbers (finitely many or infinitely many) such that for almost all $\z\in\T$,
the numbers \rf{sigma} are all nonzero distinct singular values of $\Phi(\z)$. It suffices to
prove that if $\s_{d-1}>0$, then $\Phi$ admits a partial canonical factorization of order $d$.
We prove it by induction on $d$.

Suppose first that $d=1$. Let $r=\dim{\frak S}_\Phi^{(\s_0)}(\z)$ for almost all $\z\in\T$.
Obviously,
$\dim{\frak S}_{\Phi^{\rm t}}^{(\s_0)}(\z)=r$ for almost all $\z\in\T$.
Let us show that $\Phi$ admits a factorization of of the form \rf{ordre1} with $r$-balanced
functions $\V$ and $\W^{\rm t}$. It is easy to verify that a function
$\xi\in H^2(\ell^2)$ is a maximizing vector of $H_\Phi$ if and only if
$\eta\df\bar z\ov{H_\Phi\xi}$ is a maximizing vector of $H_{\Phi^{\rm t}}$
(see \cite{P}, Ch. 14, \S 2).
Let $\M$ be the minimal invariant subspace of multiplication by $z$ on $H^2(\ell^2)$
that contains all maximizing vectors of $H_\Phi$ and 
let $\N$ be the minimal invariant subspace of multiplication 
by $z$ on $H^2(\ell^2)$ that contains all maximizing vectors of $H_{\Phi^{\rm t}}$. 

By Theorem \ref{mv}, there exist inner and co-outer functions $\U$ and $\O$ in
\lb$H^\be(\B(\C^r,\ell^2))$ such that $\M=\U H^2(\C^r)$ and $\N=\O H^2(\C^r)$.
By Theorem \ref{bal}, there exist $r$-balanced matrix functions $\V$ and $\W^{\rm t}$ 
of the form 
$$
\V=\left(\begin{array}{cc}\U&\ov{\Theta}\end{array}\right)\quad\mbox{and}\quad
\W^{\rm t}=\left(\begin{array}{cc}\O&\ov{\Xi}\end{array}\right).
$$
In exactly the same way as in the proof of Theorem 3.2 of \cite{PT3} it can be shown that
$\Phi$ admits a factorization 
$$
\Phi=\W^*\left(\begin{array}{cc}\s_0 U&\bs{0}\\\bs{0}&\Psi\end{array}\right)\V^*,
$$
where $U$ is an $r\times r$ unitary-valued matrix function. The proof of the fact that
the shift-invariant subspace spanned by the maximizing vectors of $H_U$ is $H^2(\C^r)$
is the same as it was done in the proof of Theorem 4.1 of \cite{PT3}.

In exactly the same way as in the proof of Theorem 4.1 of \cite{PT3} one can prove that
$\Psi$ satisfy condition (C). Clearly, for almost all $\z\in\T$,
$$
\s_1>\s_2>\cdots
$$
are all nonzero distinct singular values of $\Psi(\z)$.

Suppose now that $d>1$. By the inductive hypothesis, $\Psi$ admits a partial canonical
factorization of order $d-1$. Thus $\Phi$ admits a partial canonical
factorization of order $d$. $\bl$

\medskip

{\bf Remarks on uniqueness.} As we have mentioned above, unlike the 
case of finite matrix functions, in the infinite-dimensional case  a very badly approximable
function satisfying condition (C) can have infinitely many superoptimal approximants.
 However, in certain important
cases the zero function is the only superoptimal 
approximant of a very badly approximable function $\Phi$:
\begin{enumerate}
\item[(i)] if the Hankel operator $H_\Phi$ is compact, then 
$\Phi$ has a unique superoptimal approximant
(\cite{T}, \cite{Pe1}, \cite{PT1}), and so in this case 
(note that such functions $\Phi$ are automatically admissible) $\Phi$ is very badly
approximable if and only if condition (C) holds and in this case the zero function is the only
superoptimal approximant of $\Phi$;
\item[(ii)] if $\rank \Phi(\z)$ is uniformly bounded for almost all $\z\in\T$ (this happens,
for example, if $\Phi$ has finitely many columns or rows), then the family of subspaces 
${\frak S}_\Phi^{(\s)}$ stabilizes and we have the situation similar to the case 
of finite matrix functions;
in this case again the zero function is the only superoptimal approximant of $\Phi$ provided
$\Phi$ satisfies condition (C); 
\item[(iii)] if $\Phi$ satisfies condition (C) and if for almost all $\z\in\T$ the subspaces 
${\frak S}_\Phi^{(\s)}(\z)$, $\s>t_\infty(\Phi)$ span $\ell^2$, it is not hard to see that
the zero function is the only superoptimal approximant of $\Phi$. 
\end{enumerate}

\medskip

Let us explain (iii) in more detail. 

Suppose that $\Phi$ satisfies condition (C) and let $\s_k$, $k\ge0$,
be the decreasing sequence such that for almost all $\z\in\T$, the $\s_k$ are all
nonzero distinct singular values of $\Phi(\z)$
(we have already mentioned above that (C) implies that the singular values of
$\Phi(\z)$ are constant for almost all $\z\in\T$). 

It was shown in \S 5 of \cite{PT3} that if $F$ is a superoptimal approximation of $\Phi$, then
$$
\Phi(\z)-F(\z) \big|  {\frak S}_\Phi^{(\s_k)}(\z) = \Phi(\z) \big|  {\frak
S}_\Phi^{(\s_k)}(\z)\quad\mbox{for almost all}\quad\z\in\T.
$$
(this was done in \cite{PT3} for finite matrix functions, 
but the same proof also works in the infinite-dimensional case). Thus if we assume that the
subspaces ${\frak S}_\Phi^{(\s_k)}(\z)$ span $\ell^2$ for almost all
$\z\in \T$, we obtain $\Phi-F=\Phi$, and so the zero function is the only superoptimal
approximant of $\Phi$.

\

\section{\bf Badly approximable operator functions}
\setcounter{equation}{0}

\

In \cite{PT3} we obtained a description of badly approximable matrix functions.
Now we can obtain the same result for operator function.

\begin{thm}
\label{ba}
Let $\Phi\in L^\be(\B(\ell^2))$ and $\|H_\Phi\|_{\rm e}<\|\Phi\|_{L^\be}$.
If $\Phi$ is badly approximable, then
\item{\rm(i)}
$\|\Phi(\z)\|_{\B(\ell^2)}$ is constant for almost all $\z\in\T$;
\item{\rm(ii)} 
there exists a function $f$ in $\Ker T_\Phi$ such that $f(\z)$ is
a maximizing vector of $\Phi(\z)$ for almost all $\z\in\T$.

Conversely, if $\Phi\in L^\be(\B(\ell^2))$ and satisfies {\em(i)} and {\em(ii)}, then 
$\Phi$ is badly approximable. 
\end{thm}

The proof is exactly the same as in the proof of Theorem 6.1 of \cite{PT3}.

Another result of \S 6 of \cite{PT3} is a characterization of 
the set of badly approximable functions $\Phi$ such that $\|H_\Phi\|_{\rm e}<\|\Phi\|_{L^\be}$
and $\0$ is the only best approximant of $\Phi$. We can ask the same question in the case of
infinite matrix functions. However, if $\Phi\in L^\be(\B(\ell^2))$ and 
$\|H_\Phi\|_{\rm e}<\|\Phi\|_{L^\be}$, then $\0$ cannot be the only best approximant.
Indeed $\Phi$ is a badly approximable function satisfying 
$\|H_\Phi\|_{\rm e}<\|\Phi\|_{L^\be}$, then  by Theorem \ref{plokho}, it admits a partial
canonical factorization
$$
\Phi=\W^*\left(\begin{array}{cc}\s U&\0\\\0&\Psi
\end{array}\right)\V^*
$$
and $\|H_\Psi\|<\s$. Then there are infinitely many functions $Q$ in $H^\be(\B(\ell^2))$
such that $\|\Psi-Q\|_{L^\be}<\s$. Now it is easy to verify (see Theorem 1.8 of Ch. 14 of
\cite{P}) that 
$$
\Phi-\Xi Q\Theta^{\rm t}
=\W^*\left(\begin{array}{cc}\s U&\0\\\0&\Psi-Q
\end{array}\right)\V^*,
$$
where $\Theta$ and $\Xi$ are as in \rf{VW}. Thus $\Phi$ has infinitely many best approximants.

\medskip

However, we still can obtain a sufficient condition for a badly
approximable operator function to have a unique best approximant. Clearly, such a function
$\Phi$ cannot satisfy the inequality $\|H_\Phi\|_{\rm e}<\|\Phi\|_{L^\be}$.
It is convenient to normalize $\Phi$ with the condition $\|\Phi\|_{L^\be}=1$.

\begin{thm}
\label{t7.2}
Let $\Phi\in L^\infty(\B(\ell^2))$ be a function such that $\Phi(\z)$ 
is an isometry for almost all $\z\in\T$ or $\Phi^*(\z)$ is an isometry for almost all
$\z\in\T$. Suppose that $\{f(\z):f\in \Ker T_\Phi\}$ is a dense subset of  
$\Ker \Phi(\z)^\perp$ for almost all $\z\in \T$. Then  $\Phi$ is very badly approximable and
the zero function is the only best approximant of $\Phi$. 
\end{thm}

\Pf
Clearly, the fact that $\Phi$ is very badly approximable is an immediate consequence
of Theorem \ref{ba}.
Let $F$ be a best approximant of $\Phi$ and let $\Psi=\Phi-F$. 
Take $f\in \Ker T_\Phi$. 
Suppose that $\Phi(\z)$ be a coisometry for almost all $\z\in\T$. By the assumption of the
theorem $f(\z)\in \Ker \Phi(\z)^\perp$ for almost all $\z\in \T$. Hence, 
\begin{equation}
\label{7.1}
\|\Phi(\z)f(\z)\|=\|f(\z)\|\quad \text{ for almost all }\quad\z\in\T 
\end{equation}
(in the case when $\Phi(\z)$ is an isometry the above identity 
holds automatically). Since $f\in \Ker T_\Phi$, we conclude that $H_\Phi f = \Phi f$, and
\eqref{7.1} implies that 
$$
\|H_\Phi f\|_2 = \|\Phi f\|_2 = \|f\|_2. 
$$
Since $F\in H^\infty$, we have $H_\Phi=H_\Psi$. Consider the
following chain of inequalities:
$$
\|f\|_2=\|H_\Phi f\|_2 =\|H_\Psi f\|_2\le \|H_\Psi f\|_2\le \|\Psi\|_\infty\|f\|_2 \le \|f\|_2.
$$
Therefore all inequalities in this chain are, in fact, equalities, and so
$$
\Psi f = H_\Psi f = H_\Phi f = \Phi f.
$$
Since the set $\{f(\z):f\in \Ker T_\Phi\}$ is a dense subset of  
$\Ker \Phi(\z)^\perp$ for almost all $\z\in \T$, we obtain 
\begin{equation}
\label{7.2}
\Phi(\z)\big|\Ker\Phi(\z)^\perp=\Psi(\z)\big|\Ker\Phi(\z)^\perp\quad\text{for almost
all }\quad\z\in \T. 
\end{equation}

If $\Phi(\z)$ is an isometry for almost all $\z\in\T$, then $\Ker \Phi(\z)$ is trivial, and 
therefore $\Phi=\Psi$. 

If $\Phi(\z)$ is a coisometry for almost all $\z\in\T$, then \eqref{7.2} implies that for
almost all $\z\in \T$
$$
(\Phi(\z)x, y) = (\Psi(\z)x, y),\quad x\in \Ker \Phi(\z)^\perp,~y\in \ell^2. 
$$
Since  $\|\Psi\|_\infty\le 1$, it follows that $\Phi^*(\z)=\Psi^*(\z)$ 
for almost all $\z\in\T$. 
$\bl$

\

\

\noindent
\begin{tabular}{p{8cm}p{14cm}}
V.V. Peller & S.R. Treil \\
Department of Mathematics & Department of Mathematics \\
Michigan State University  & Brown University \\
East Lansing, Michigan 48824 & Providence, Rhode Island 02912\\
USA&USA
\end{tabular}

\end{document}